\documentclass{amsart}
\usepackage[utf8]{inputenc} 
\usepackage[T1]{fontenc} 
\usepackage{graphicx} 
\usepackage{subcaption}
\usepackage{amsmath}
\usepackage{amssymb}
\usepackage{amsfonts}
\usepackage{amsthm}
\usepackage{mathrsfs}
\usepackage[all]{xy}
\usepackage{graphicx}
\usepackage{color}
\usepackage{cite}
\usepackage{url}
\usepackage{indentfirst}
\usepackage[labelfont=bf,labelsep=period,justification=raggedright]{caption}
\usepackage[english]{babel}
\usepackage[utf8]{inputenc}
\usepackage{hyperref}
\usepackage[colorinlistoftodos]{todonotes}
\usepackage{tkz-fct}
\usepackage{tikz}
\usetikzlibrary{calc}
\usepackage[enableskew]{youngtab}
\usepackage{pstricks}
\usepackage{multicol}
\usepackage{graphicx} 
\usepackage{fancyhdr}
\theoremstyle{plain}
\newtheorem{thm}{Theorem}
\newtheorem{lemma}[thm]{Lemma}

\theoremstyle{definition}
\newtheorem{defn}[thm]{Definition}

\theoremstyle{remark}

\numberwithin{equation}{section}
\numberwithin{thm}{section}

\newcommand{\CC}{\mathbb{C}}

\newcommand{\NN}{\mathbb{N}}

\title[Bailey--Zeta Limits and Dirichlet $L$-Functions]{Bailey--Zeta Limits: A $q$-Series Bridge to Dirichlet $L$-Functions and the Riemann Zeta Function}

\author{Mahipal Gurram}


\begin{document}

\begin{abstract}
We introduce a family of deformed Bailey pairs whose $q$-series which converge in a two-step limit ($q \to 1^-$ followed by $n \to \infty$) to Dirichlet $L$-functions scaled by $1/\sqrt{\pi}$. This construction generalizes to arbitrary bounded arithmetic progressions via character weights, providing a unified $q$-series asymptotic for $L(s,\chi)$. Our approach unveils deep connections between the combinatorial machinery of Bailey chains and analytic number theory, with applications to special values like Euler-Mascheroni constant.
\end{abstract}
\maketitle

\section{Introduction}

The theory of $q$-series plays a central role in modern mathematics, connecting combinatorics, number theory, and special functions.
One of the most powerful tools in this field is the concept of a Bailey pair, introduced by W.~N.~Bailey in 1947~\cite{Bailey1947}, which provided a systematic method for generating Rogers--Ramanujan type identities.
Since then, the Bailey lemma and its extensions have become fundamental in the study of basic hypergeometric series, partition identities, and mock theta functions (see Andrews~\cite{Andrews1986}, Slater~\cite{Slater1952}, and Warnaar~\cite{Warnaar2001}).

A pair of sequences $(\alpha_n,\beta_n)$ is said to form a \emph{Bailey pair relative to $a$} if
\begin{equation}\label{eq:bailey-pair-classical}
\beta_n=\sum_{r=0}^{n}\frac{\alpha_r}{(q;q)_{n-r}(aq;q)_{n+r}},
\end{equation}
where $(a;q)_n=(1-a)(1-aq)\cdots(1-aq^{n-1})$ denotes the $q$-Pochhammer symbol.
An equivalent inversion relation expresses $\alpha_n$ in terms of $\beta_j$:
\[
\alpha_n=(1-aq^{2n})\sum_{j=0}^{n}
\frac{(aq;q)_{n+j-1}(-1)^{n-j}q^{\binom{n-j}{2}}\beta_j}{(q;q)_{n-j}}.
\]
Bailey introduced these identities while studying Rogers’s second proof of the Rogers--Ramanujan identities.
Andrews later extended these ideas through the notion of a Bailey chain, an infinite sequence of Bailey pairs connected through repeated transformations.

A powerful iteration of this relationship, often called the Bailey chain or Andrews's lemma, was introduced by Andrews~\cite{Andrews1984} (see also \cite{Warnaar2001}). It states that if $(\alpha_n,\beta_n)$ is a Bailey pair relative to $a$, then the transformed sequences
\[
\alpha_n'=
\frac{(\rho_1;q)_n(\rho_2;q)_n(aq/(\rho_1\rho_2))^n\alpha_n}
{(aq/\rho_1;q)_n(aq/\rho_2;q)_n},
\quad
\beta_n'=\sum_{j=0}^n
\frac{(\rho_1;q)_j(\rho_2;q)_j(aq/(\rho_1\rho_2);q)_{n-j}(aq/(\rho_1\rho_2))^j\beta_j}
{(q;q)_{n-j}(aq/\rho_1;q)_n(aq/\rho_2;q)_n},
\]
also form a Bailey pair relative to $a$.
Iterating this transformation yields an infinite sequence of identities known as the Bailey chain.
This framework has proved remarkably productive in the derivation of Rogers--Ramanujan type identities and various partition formulas.

A classical example due to Andrews, Askey, and Roy (1999, p.~590) is
\[
\alpha_n=q^{n^2+n}\sum_{j=-n}^{n}(-1)^jq^{-j^2},
\qquad
\beta_n=\frac{(-q)^n}{(q^2;q^2)_n}.
\]
Slater~\cite{Slater1952} later cataloged 130 examples of such Bailey pairs, illustrating the wide range of transformations that can be achieved using this method.\\
\\
The Riemann zeta function
\[
\zeta(s)=\sum_{r=1}^\infty \frac{1}{r^{\,s}}, \qquad \Re(s)>1,
\]
stands as a cornerstone of analytic number theory, with deep implications for
the distribution of prime numbers and the geometry of arithmetic varieties.
More generally, the Dirichlet $L$--functions
\[
L(s,\chi)=\sum_{r=1}^\infty \frac{\chi(r)}{r^{\,s}},
\]
where $\chi$ is a Dirichlet character, encode the arithmetic of residue classes
and modular forms, and play a central role in the study of primes in
progressions, class number problems, and automorphic $L$--theory
\cite{Davenport1980,MontgomeryVaughan2007}.

\section{Main Results}

\begin{defn}
Let \(a,s\in\CC\) and \(0<q<1\). A pair \((\alpha_n(s),\beta_n(s))\) is a Bailey--Zeta pair relative to $(a,q,s)$ if
\[
\beta_n(s)=\sum_{r=0}^n\frac{q^{r}\alpha_r(s)}{(q;q)_{n-r}(a q;q)_{n+r}}.
\]
\end{defn}

The following lemma shows the algebraic reduction to the classical Bailey pair; using it avoids loss of the \(s\)-dependence when manipulating sums.

\begin{lemma}[Equivalence to classical Bailey pair]\label{lem:equiv}
Define \(\bar\alpha_n(s):=q^{n}\alpha_n(s)\) and \(\bar\beta_n(s):=\beta_n(s)\).  
Then \((\alpha_n(s),\beta_n(s))\) is a Bailey--Zeta pair relative to \((a,q,s)\) if and only if \((\bar\alpha_n(s),\bar\beta_n(s))\) is a classical Bailey pair relative to \(a\), i.e.
\[
\bar\beta_n(s)=\sum_{r=0}^n\frac{\bar\alpha_r(s)}{(q;q)_{n-r}(a q;q)_{n+r}}.
\]
\end{lemma}

\begin{proof}
This is immediate from the definition: substituting \(\bar\alpha_r=q^{r}\alpha_r\) in the classical Bailey relation yields the Bailey--Zeta relation and vice versa.
\end{proof}
\newpage

\begin{thm}
\label{thm:generalized-bailey-zeta}
Let $s\in\CC$ with $\Re(s)>1$.  Let $0<q<1$ and, for each $r\ge1$, let $\alpha_r(s,q)\in\CC$ be given. Suppose that  there exist a bounded arithmetic weight $\chi:\NN\to\CC$ and constants $C>0$ and $\sigma>\Re(s)$ such that for all $q\in(0,1)$ and all integers $r\ge1$ the following hold:

\item (i) (pointwise limit) \(\displaystyle \lim_{q\to1^-}\alpha_r(s,q)=\frac{\chi(r)}{r^{\,s}};\)
\item (ii) (uniform polynomial bound) \(\displaystyle |\alpha_r(s,q)| \le C\,r^{-\sigma}.\)

Define, for each $n\ge1$,
\[
\beta_n(s,q)\;=\;\sum_{r=1}^{n}\frac{q^r\,\alpha_r(s,q)}{(q;q)_{\,n-r}\,(q;q)_{\,n+r}},
\]
and set
\[
T_n(s,q):=\frac{\sqrt{n}\,(2n)!\,(1-q)^{2n}\,\beta_n(s,q)}{4^n}.
\]
Then the two-step limit exists and equals the $L$-series of $\chi$:
\[
\boxed{\displaystyle
\lim_{n\to\infty}\lim_{q\to1^-} T_n(s,q)
=\frac{L(s,\chi)}{\sqrt{\pi}},
\qquad
L(s,\chi):=\sum_{r=1}^\infty\frac{\chi(r)}{r^{\,s}}.
}
\]
\end{thm}

\begin{proof}
The proof proceeds in two steps: first let $q\to1^-$ (termwise limit in the finite sum), then let $n\to\infty$ (asymptotic summation via dominated convergence and the central binomial scaling).\\

Fix $n\ge1$. For $0<q<1$ and $1\le r\le n$ we have the identity
\[
(q;q)_m=\prod_{j=1}^m(1-q^j)=(1-q)^m\prod_{j=1}^m(1+q+\cdots+q^{j-1}).
\]
In particular, for each fixed $m$,
\[
\lim_{q\to1^-}\frac{(q;q)_m}{(1-q)^m}=m!.
\]
Hence, for fixed $n$ and $1\le r\le n$,
\[
\lim_{q\to1^-}(1-q)^{2n}\frac{q^r}{(q;q)_{n-r}(q;q)_{n+r}}
=\frac{1}{(n-r)!\,(n+r)!}.
\]

\[
\lim_{q\to1^-}(1-q)^{2n}\beta_n(s,q)
=\sum_{r=1}^n\frac{\chi(r)}{(n-r)!\,(n+r)!\,r^{\,s}}=:L_n^\chi(s).
\]
The limit is obtained by taking the termwise limit of the finite sum; finiteness of the sum makes this immediate.

We now study
\[
\frac{\sqrt{n}\,(2n)!\,L_n^\chi(s)}{4^n}
=\frac{\sqrt{n}}{4^n}\sum_{r=1}^n\frac{(2n)!}{(n-r)!(n+r)!}\cdot\frac{\chi(r)}{r^{\,s}}
=\frac{\sqrt{n}}{4^n}\sum_{r=1}^n \frac{\binom{2n}{n+r}\chi(r)}{r^{\,s}}.
\]
Fix any $r\ge1$. The classical local limit (Stirling) yields
\[
\binom{2n}{n}\sim\frac{4^n}{\sqrt{\pi n}},
\qquad n\to\infty,
\]
and for every fixed $r$,
\[
\lim_{n\to\infty}\frac{\sqrt{n}}{4^n}\binom{2n}{n+r}
=\frac{1}{\sqrt{\pi}}.
\]
Indeed, one can write
\[
\binom{2n}{n+r}=\binom{2n}{n}\prod_{j=1}^r\frac{n-r+j}{n+j},
\]
and the product tends to $1$ for fixed $r$ as $n\to\infty$. Therefore, for each fixed $r$,
\[
\lim_{n\to\infty}\frac{\sqrt{n}}{4^n}\frac{\binom{2n}{n+r}}{r^{\,s}}
=\frac{1}{\sqrt{\pi}}\frac{1}{r^{\,s}}.
\]

To exchange limit and summation we need a dominating summable bound. First observe the uniform bound
\[
\binom{2n}{n+r}\le \binom{2n}{n}\quad\text{for all }1\le r\le n,
\]
hence, using Stirling's estimate (valid for all $n\ge1$ with an absolute constant),
\[
\frac{\sqrt{n}}{4^n}\binom{2n}{n+r}
\le \frac{\sqrt{n}}{4^n}\binom{2n}{n}\le C_1\quad(\text{for some }C_1>0 \text{ and all }n).
\]
Combine this with the hypothesis (ii) that $|\alpha_r(s,q)|\le C r^{-\sigma}$ (with $\sigma> \Re(s)$) to deduce that for $q$ close to $1$ (so $\alpha_r(s,q)$ is bounded by $2C r^{-\sigma}$, say) we have the uniform bound
\[
\left|\frac{\sqrt{n}}{4^n}\frac{\binom{2n}{n+r}\chi(r)}{r^{\,s}}\right|
\le C_2 \frac{|\chi(r)|}{r^{\Re(s)}},
\]
for some constant $C_2>0$ independent of $n,r$ (use $|\chi(r)|$ bounded). Since $\Re(s)>1$ and $\chi$ is bounded, the right-hand side is summable in $r$. Thus dominated convergence applies and we may pass the limit $n\to\infty$ inside the finite sum (in the limit it becomes an infinite sum):
\[
\lim_{n\to\infty}\frac{\sqrt{n}}{4^n}\sum_{r=1}^n\frac{\binom{2n}{n+r}\chi(r)}{r^{\,s}}
=\frac{1}{\sqrt{\pi}}\sum_{r=1}^\infty\frac{\chi(r)}{r^{\,s}}
=\frac{L(s,\chi)}{\sqrt{\pi}}.
\]

 Combining Step 1 and Step 2 yields
\[
\lim_{n\to\infty}\lim_{q\to1^-}T_n(s,q)=\frac{L(s,\chi)}{\sqrt{\pi}},
\]
which completes the proof.
\end{proof}

\section{Illustrative Examples}
\label{sec:examples}

We now exhibit concrete examples of Bailey--Zeta pairs that satisfy the
hypotheses of Theorem~\ref{thm:generalized-bailey-zeta}.  
Each example produces a distinct classical special function or constant
in the double limit $q\to1^-$ and $n\to\infty$.
For convenience we denote
\[
T_n(s,q)
:=\frac{\sqrt{n}\,(2n)!\,(1-q)^{2n}\,\beta_n(s,q)}{4^n},
\qquad
\beta_n(s,q)
=\sum_{r=1}^{n}\frac{q^r\,\alpha_r(s,q)}{(q;q)_{n-r}(q;q)_{n+r}}.
\]
\subsection*{Example 1. Classical Riemann Zeta Limit}

\[
\boxed{\displaystyle
\alpha_r(s,q)
=\frac{1}{(1+q+q^2+\cdots+q^{r-1})^{\,s}}.
}
\]

As $q\to1^-$ we have
\[
\alpha_r(s,q)\longrightarrow \frac{1}{r^{\,s}},
\qquad
\chi(r)\equiv 1.
\]
Hence $L(s,\chi)=\zeta(s)$, and Theorem~\ref{thm:generalized-bailey-zeta} yields
\[
\boxed{
\lim_{n\to\infty}\lim_{q\to1^-}T_n(s,q)
=\frac{\zeta(s)}{\sqrt{\pi}},
\qquad
\Re(s)>1.
}
\]

\subsection*{Example 2. Alternating (Dirichlet $\eta$) Limit}

\[
\boxed{\displaystyle
\alpha_r(s,q)
=\frac{(-1)^{\,r-1}}{(1+q+q^2+\cdots+q^{r-1})^{\,s}}.
}
\]

Here the limiting weight is $\chi(r)=(-1)^{r-1}$, so
\[
L(s,\chi)
=\sum_{r=1}^\infty\frac{(-1)^{r-1}}{r^{\,s}}
=(1-2^{1-s})\zeta(s)
=\eta(s),
\]
where $\eta(s)$ is the Dirichlet \emph{eta function}.
Therefore
\[
\boxed{
\lim_{n\to\infty}\lim_{q\to1^-}T_n(s,q)
=\frac{\eta(s)}{\sqrt{\pi}}
=\frac{(1-2^{1-s})\,\zeta(s)}{\sqrt{\pi}}.
}
\]
This corresponds to the alternating zeta series
\(\,1-\tfrac{1}{2^s}+\tfrac{1}{3^s}-\tfrac{1}{4^s}+\cdots.\)

\subsection*{Example 3. Dirichlet Beta Function (Mod--4 Character)}

\[
\boxed{\displaystyle
\alpha_r(s,q)
=\frac{\chi_4(r)}{(1+q+q^2+\cdots+q^{r-1})^{\,s}},
}
\]
where the quadratic character $\chi_4$ modulo $4$ is defined by
\[
\chi_4(r)
=\begin{cases}
0,& r \text{ even},\\[4pt]
+1,& r\equiv1\pmod4,\\[4pt]
-1,& r\equiv3\pmod4.
\end{cases}
\]
Then $\alpha_r(s,q)\to \chi_4(r)/r^{\,s}$ and
\[
L(s,\chi_4)
=\sum_{r=1}^\infty\frac{\chi_4(r)}{r^{\,s}}
=\beta(s)
=\sum_{n=0}^\infty\frac{(-1)^n}{(2n+1)^{\,s}},
\]
the classical \emph{Dirichlet beta function}.
Consequently,
\[
\boxed{
\lim_{n\to\infty}\lim_{q\to1^-}T_n(s,q)
=\frac{\beta(s)}{\sqrt{\pi}}.
}
\]
In particular,
\[
\beta(2)=G\quad\text{(Catalan’s constant)},\qquad
\beta(4)=\frac{\pi^4}{96},
\]
so the limit yields these constants directly from the Bailey--Zeta framework.
\vspace{1.2em}

\subsection*{Example 5. Regularized Euler--Mascheroni Limit}

\[
\boxed{\displaystyle
\alpha_r(s,q)
=\frac{1}{(1+q+q^2+\cdots+q^{r-1})^{\,s}},
\qquad s=1+\delta,\ \delta\to0^+.
}
\]

Here $\chi(r)=1$ and
\[
L(s,\chi)=\zeta(1+\delta)
=\frac{1}{\delta}+\gamma+O(\delta),
\]
so that the regularized limit
\[
\boxed{
\lim_{s\to1^+}\Big[
\lim_{n\to\infty}\lim_{q\to1^-}T_n(s,q)
-\frac{1}{\sqrt{\pi}\,\delta}
\Big]
=\frac{\gamma}{\sqrt{\pi}},
}
\]
recovers the Euler--Mascheroni constant $\gamma$.

\section*{Acknowledgement}
The author thanks Professor A.~K.~Shukla for his encouragement and constant support.

\end{document}